\documentclass[10pt]{article}
\usepackage{amssymb,amsmath,mathrsfs,latexsym,amscd}
\usepackage{exscale,latexsym,amsthm,graphics}

\usepackage{makeidx}
\makeindex
\usepackage{epsfig}
\usepackage{txfonts}

\newtheorem{thm}{Theorem}[section]

\newtheorem{lemma}[thm]{Lemma}
\newtheorem{prop}[thm]{Proposition}
\newtheorem{proposition}[thm]{Proposition}
\newtheorem{defn}[thm]{Definition}

\newtheorem{remark}[thm]{Remark}
\newtheorem{example}[thm]{Example}

\def\bM {{\mathbb M}}  
 \def\Q {{\mathbb Q}}

\def\N {{\mathbb N}}

\def\qed{{\hfill $\Box$ \bigskip}}
\def\N {{\mathbb N}}
\def\R {{\mathbb R}}

\def\EE{{\mathbb E}}
\def\P{{\mathbb P}}
\newcommand{\F}{\mathcal{F}}

\def\E{{\mathcal E}}
\def\P{{\mathbb P}}

\begin{document}
\noindent
{{\Large\bf Pointwise weak existence for diffusions associated with degenerate elliptic forms and 2-admissible weights}\footnote{This research was supported by NRF-DFG Collaborative Research program and Basic Science Research Program 
through the National Research Foundation of Korea (NRF-2012K2A5A6047864 and NRF-2012R1A1A2006987) 
and by DFG through Grant Ro 1195/10-1.}}

\bigskip
\noindent
{\bf Jiyong Shin},
{\bf Gerald Trutnau}
\\

\noindent
{\small{\bf Abstract.} 
Using analysis for 2-admissible functions in weighted Sobolev spaces and stochastic calculus for possibly degenerate symmetric elliptic forms, we construct weak solutions to a wide class of stochastic differential equations starting from an explicitly specified subset in Euclidean space. The solutions have typically unbounded and discontinuous drift but may still in some cases start from all points of $\R^d$ and thus in particular from those where the drift terms are infinite. As a consequence of our approach we are able to provide new non-explosion criteria for the unique strong solutions of \cite{Zh}.  \\



\section{Introduction}
In this paper, we consider a symmetric Dirichlet form (given as the closure of)
\[
\E^A(f,g): = \frac{1}{2} \int_{\R^d} \langle A  \nabla f , \nabla g \rangle \, dx, \quad f,g \in C_0^{\infty}(\R^d)
\]
on $L^2(\R^d,m)$, $m:= \rho dx$. Such forms were considered under analytic aspects in \cite{FJK}, \cite{FKS} and \cite{BM}. The basic conditions on $\rho$ and on the diffusion matrix $A = (a_{ij})_{1 \le i,j \le d}$ are formulated in (I)-(IV) and (HP1)-(HP2) below. 

Our first aim is to construct an associated Hunt process to $\E^A$ which satisfies Fukushima's absolute continuity condition 
(cf. Remark \ref{ch5;remsuphu} (i) and for the absolute continuity condition \cite[(4.2.9) and Theorem 5.5.5]{FOT} and \cite{fuku94}) 
and subsequently to identify the stochastic differential equation (hereafter SDE) with explicit form \eqref{ch5;sfd2} verified by it for as much as possible explicitly given starting points $x \in \R^d$. This is done under some additional assumptions, namely (HP3), \eqref{ch5;emucpol} or (HP3)$^{\prime}$, \eqref{ch5;eqmupo2} in Section \ref{ch5;sect3}  and (HP4), (HP5) in Section \ref{ch5;sect4}. Here we follow the major lines of the program developed in \cite{ShTr13a}, which explicitly provides tools to apply Fukushima's absolute continuity condition. 
In the situation of \cite[Section 3]{ShTr13a} it is a well-known fact that the intrinsic metric (from the Dirichlet form there) is equal to the Euclidean metric.
However, in contrast to \cite[Section 3]{ShTr13a} the handling of more general 2-admissible weights $\rho$ 
(including Muckenhoupt $A_2$-weights appearing in \cite[Section 3]{ShTr13a}) 
and a possibly degenerate diffusion matrix $A$ requires (strong) equivalence between the intrinsic metric (derived from $\E^A$) and the Euclidean metric 
(see Lemmas \ref{ch5;l;intrinsic}, \ref{l;ch5;doublepseu}, \ref{ch5;l;poincom}, \ref{ch5;l;poincare} and Remark \ref{ch5;rebmea}).
 These generalized assumptions on $\rho$ and $A$ then extend results from \cite{ShTr13a} (cf. Remark \ref{ch5;reexshtr}).

Our second aim is to provide new non-explosion results (cf. Remark \ref{ch5;renonexpn}) for the solution to the SDE
\[
X_t = x + \int_0^t \sigma(X_s) \, dW_s + \int_0^t b(X_s) \, ds, \quad x \in \R^d,
\]
where $\sigma$ satisfies (C1)-(C3) of \cite{Zh} and $b \in L^{2(d+1)}_{loc}(\R^d,dx)$. In particular, applying \cite[Theorem 1.1]{Zh} to our solution of \eqref{ch5;sfd2} and using Dirichlet form theory, we obtain under the conditions stated in Theorem \ref{ch5;t;ssoleae} that \eqref{ch5;sfd2} admits a unique strong solution which is non-explosive. This completes our results from \cite{RoShTr} where we presented new non-explosion results for the unique strong solutions of \cite{KR} (see also \cite{FF}). Finally, let us mention that the results of \cite{BGS} and \cite{RoShTr} are particularly close and complementary to ours. Moreover, our results can also be used to obtain new dynamics for interacting particle systems (cf. \cite{AKR} and \cite{BGS}).

The organization of the paper is as follows. Right after the introduction in Section \ref{ch5;2;prelide}, we develop the framework and analytic background based on results from \cite{BM,  Ca,  FJK, FKS, GrHu, HK, HKM, ShTr13a, ShTr13b, St3, St5}. In Section \ref{ch5;sect3} we construct a Hunt process satisfying the absolute continuity condition and identify it with weak solutions related to concrete 2-admissible weights of the form \eqref{ch5;eq;2admie} below. In Section \ref{ch5;sect4} we do the same for weights $\rho$ in a subclass of the Muckenhoupt $A_2$-class satisfying (HP4). Section \ref{ch5;sect5} is devoted to the new non-explosion results.

\section{Preliminaries and degenerate elliptic forms with respect to 2-admissible weights}\label{ch5;2;prelide}

For $E \subset \R^d$, $d \ge 2$ with Borel  $\sigma$-algebra $\mathcal{B}(E)$ we denote the set of all $\mathcal{B}(E)$-measurable $f : E \rightarrow \R$ which are bounded by $\mathcal{B}_b(E)$. The usual $L^q$-spaces $L^q(E, \mu)$, $q \in[1,\infty]$ are equipped with $L^{q}$-norm $\| \cdot \|_{q}$ with respect to the  measure $\mu$ on $E$, $\mathcal{A}_b$ : = $\mathcal{A} \cap \mathcal{B}_b(E)$ for $\mathcal{A} \subset L^q(E,\mu)$,  and $L^{q}_{loc}(E,\mu) := \{ f \,|\; f \cdot 1_{U} \in L^q(E, \mu),\,\forall U \subset E, U \text{ relatively compact open} \}$, where $1_A$ denotes the indicator function of a set $A$. If $\mathcal{A}$ is a set of functions $f : E \to \R$, we define $\mathcal{A}_0 : = \{f \in \mathcal{A} \ | $ supp($f$) : = supp($|f| m$) is compact in $E \}$. Let $\nabla f : = ( \partial_{1} f, \dots , \partial_{d} f )$ where $\partial_j f$ is the $j$-th weak partial derivative of $f$ and $\partial_{ij} f := \partial_{i}(\partial_{j} f) $, $i,j=1, \dots, d$. For any open set $E \subset \R^d$ we also denote the set of continuous functions on $E$, the set of continuous bounded functions on $E$, the set of compactly supported continuous functions in $E$ by $C(E)$, $C_b(E)$, $C_0(E)$, respectively.  For any open set $E \subset \R^d$ $C_{\infty}(E)$ denotes the space of continuous functions on $E$ which vanish at infinity and $C^{\infty}(E)$, $C_0^{\infty}(E)$ denote the set of all infinitely differentiable functions on $E$,  the set of all infinitely differentiable functions with compact support in $E$, respectively. As usual $dx$ denotes Lebesgue measure on $\R^d$ and for any open set $E \subset \R^d$ the Sobolev space $H^{1,q}(E, dx)$, $q \ge 1$ is defined to be the set of all functions $f \in L^{q}(E, dx)$ such that $\partial_{j} f \in L^{q}(E, dx)$, $j=1, \dots, d$, and $H^{1,q}_{loc}(\R^d, dx) : =  \{ f  \,|\;  f \cdot \varphi \in H^{1,q}(\R^d, dx),\,\forall \varphi \in  C_0^{\infty}(\R^d)\}$. We always equip $\R^d$ with the Euclidean norm $\| \cdot \|$ with corresponding inner product $\langle \cdot, \cdot \rangle$ and write $B_{r}(x): = \{ y \in \R^d \ | \ \|x-y\| < r  \}$, $x \in \R^d$, $r>0$. For $A \subset \R^d$ let $\overline{A}$ denote the closure of $A$ in $\R^d$.\\

We say that a locally integrable function $\rho: \R^d \to \R$ is $2$-admissible if the following four conditions are satisfied (see \cite[Section 1.1]{HKM}):
\begin{itemize}
\item[(I)]
$0 < \rho(x) < \infty$ for $dx$-a.e. $x \in \R^d$ and $\rho$ is doubling, i.e. there is a constant $C_1 > 0$ such that
for any   ball $B_{r}(x) \subset \R^d$
\begin{equation}\label{ch5;doublep}
m(B_{2r}(x))  \le  C_1 \  m(B_{r}(x)), \quad m : = \rho dx.
\end{equation}
\item[(II)] If $D \subset \R^d$ is an open set and $(u_n)_{n \ge 1} \subset C^{\infty}(D)$ is a sequence of functions such that
\[
\lim_{n \to \infty} \int_{D} |u_n|^2 \, dm =0, \quad \text{and} \ \lim_{n \to \infty} \int_{D} \|\nabla u_n - \vartheta \|^2 \, dm =0,
\]
where $\vartheta$ is an $\R^d$-valued measurable function with $\|\vartheta\|\in L^2(D,m)$, then $\vartheta=0$. 
\item[(III)] There are constants $\theta > 1$ and $C_2 > 0$ such that for $x \in \R^d$ and $r>0$
\[
\left( \frac{1}{m(B_r(x))} \int_{B_r(x)} |u|^{2 \theta } \, dm  \right)^{1 / {2 \theta}} \le C_2 \ r  \left( \frac{1}{m(B_r(x))} \int_{B_r(x)} \| \nabla  u \|^2 \, dm \right)^{1/2}, \ \forall u \in C_0^{\infty}(B_r(x)).
\]
\item[(IV)] There is a constant $C_3 > 0$ such that for $x \in \R^d$ and $r>0$
\[
\int_{B_r(x)} |u -u_{x,r}|^2 \, dm \le C_3 \ r^2 \ \int_{B_r(x)} \| \nabla u \|^2 \, dm, \quad \forall u \in C^{\infty}(B_r(x)) \cap C_b( B_r(x)),
\]
where $u_{x,r} = \frac{1}{m(B_r(x))}\int_{B_r(x)} u \, dm$.
\end{itemize}
\begin{remark}\label{ch5;re2admm}
We know from \cite[Theorem 13.1]{HK} that a locally integrable weight $\rho$ is a 2-admissible weight if and only if $\rho$ is doubling and there exist constants $c>0$, $\gamma \ge 1$ such that for $x \in \R^d$ and $r>0$
\[
\frac{1}{m(B_r(x))} \int_{B_r(x)} |u - u_{x,r}| \, dm \le c \ r \left( \frac{1}{m(B_{\gamma r}(x))} \int_{ B_{\gamma r}(x)} \|\nabla u \|^2 \, dm  \right)^{1/2}, 
\]
whenever $u \in C^{\infty}(B_{\gamma r}(x))$ (weak (1,2) Poincar\'{e} inequality).
\end{remark}

Throughout the whole article let $\rho$ be a locally integrable 2-admissible weight. For later use we define a symmetric bilinear form
\begin{equation}\label{ch5;dfrho}
\E^{\rho}(f,g) = \frac{1}{2} \int_{\R^d} \langle \nabla f , \nabla g \rangle \, dm, \quad f,g \in C_0^{\infty}(\R^d).
\end{equation}
By (II) $(\E^{\rho}, C_0^{\infty}(\R^d) )$ is closable  in $L^2(\R^d, m)$ and its closure $(\E^{\rho}, D(\E^{\rho}) )$ is a strongly local, regular  Dirichlet form in the sense of \cite{FOT}. 

Consider the following assumption:
\begin{itemize}
\item[(HP1)]  
Let $A=(a_{ij})_{1 \le i,j \le d}$ be a symmetric (possibly) degenerate (uniformly weighted)  elliptic $d \times d$ matrix, that is  $a_{ij} \in L^1_{loc} (\R^d,dx)$ and there exists a constant $\lambda \ge 1$
 such that for $dx$-a.e. $x \in \R^d$
\begin{equation}\label{ch5;eq;uelliptic} 
\lambda^{-1} \ \rho(x) \ \| \xi \|^2 \le \langle A(x) \xi, \xi \rangle \le \lambda \ \rho(x) \ \|\xi\|^2, \quad  \forall \xi \in \R^d.
\end{equation} 
\end{itemize}

From now on, we fix $A= (a_{ij})_{1 \le i,j \le d}$ satisfying (HP1) and consider the symmetric bilinear form
\[
\E^A(f,g) = \frac{1}{2} \int_{\R^d} \langle A \nabla f , \nabla g \rangle \, dx, \quad f,g \in C_0^{\infty}(\R^d).
\]
By closability of $(\E^{\rho}, C_0^{\infty}(\R^d) )$ in $L^2(\R^d, m)$ and \eqref{ch5;eq;uelliptic} $(\E^A,C_0^{\infty}(\R^d))$ is closable in $L^2(\R^d, m)$. The closure $(\E^A,D(\E^A))$ of $(\E^A,C_0^{\infty}(\R^d))$ is a strongly local, regular, symmetric Dirichlet form (see  \cite{FOT}). The Dirichlet form $(\E^A,D(\E^A))$ can be written as 
\[
\E^A(f,g) = \frac{1}{2} \int_{\R^d} \, d\Gamma^{A} (f,g), \quad f,g \in D(\E^A),
\]
where $\Gamma^{A}$ is a symmetric bilinear form on $D(\E^A) \times D(\E^A)$ with values in the signed Radon measures on $\R^d$ (called energy measures). By an approximation argument we can extend the quadratic form $u \mapsto \Gamma^{A}(u,u)$ to
 $ D(\E^A)_{loc} = \{ u \in L^2_{loc}(\R^d, m) \ | \ \Gamma^{A}(u,u) \ \text{is a Radon measure} \}$, where $D(\E^A)_{loc}$ is the set of all $m$-measurable functions $u$ on $\R^d$ for which on every compact set $K \subset \R^d$ there exists a function $v \in D(\E^A)$ with $u=v$ $m$-a.e on $K$ (cf. \cite[p. 274]{St3}). Then the energy measure $\Gamma^{A}$ defines in an intrinsic way a pseudo metric $d$ on $\R^d$ by
\[
d(x,y) = \sup \Big\{u(x) - u(y) \ | \ u \in D(\E^A)_{loc} \cap C(\R^d), \ \Gamma^{A}(u,u)   \le   m \ \text{on} \ \R^d  \Big\}, 
\]
where $\Gamma^{A}(u,u)   \le   m$ means that the energy measure $\Gamma^{A}(u,u)$ is absolutely continuous w.r.t. the reference measure $m$ with Radon-Nikodym derivative $\frac{d}{dm} \Gamma^{A}(u,u) \le 1$. We define the balls w.r.t. intrinsic metric by 
\[
\tilde{B}_r(x) = \{ y \in \R^d \ | \ d(x,y) < r \}, \quad x \in \R^d, \quad r>0.
\]
Let  $(T_t)_{t > 0}$ and $(G_{\alpha})_{\alpha > 0}$ be the $L^2(\R^d, m)$-semigroup and resolvent associated to $(\E^A,D(\E^A))$ (see \cite{FOT}).\\

 We assume from now on
\begin{itemize}
\item[(HP2)]  
Either $\sqrt{\rho} \in H^{1,2}_{loc}(\R^d,dx) \quad $  or $ \quad \rho^{-1} \in L^1_{loc}(\R^d,dx)$.
\end{itemize}

\begin{lemma}\label{ch5;l;intrinsic}
For any $x,y \in \R^d$
\begin{equation}\label{ch5;eq;intrinsicm}
\frac{1}{\sqrt{\lambda}} \ \| x-y \| \le d (x,y) \le \sqrt{\lambda} \ \| x-y \|   ,
\end{equation}
where $\lambda \in [1, \infty)$ is as in \eqref{ch5;eq;uelliptic}.
\end{lemma}
\proof
For any $z \in \R^d$ the map
\[
u \ : \ x \longmapsto u(x) : = \langle x, z \rangle 
\]
lies in $D({\E^A})_{loc} \cap C(\R^d)$. For fixed $y, y^{\prime} \in \R^d$ ($y \neq y^{\prime}$), choose 
\[
z = \frac{(y - y^{\prime})}{\sqrt{\lambda} \ \| y - \ y^{\prime} \|} \in \R^d.
\] 
Then by \eqref{ch5;eq;uelliptic}
\[
\int_B \Gamma^{A}(u,u)  = \int_B  \langle A \nabla u,  \nabla u \rangle dx \le \lambda \int_B   \| \nabla u\|^2 \rho \, dx = \int_B  \rho dx, \quad \forall B \in \mathcal{B}(\R^d).
\]
Hence $\Gamma^{A}(u,u)   \le  m $.
Furthermore
\[
u(y) - u(y^{\prime}) = \frac{1}{\sqrt{\lambda}} \ \| y- y^{\prime} \|.
\]
Therefore for any $x,y \in \R^d$
\[
d (x,y) \ge \frac{1}{\sqrt{\lambda}} \ \| x-y \|.
\]
Define 
\[
d^{\rho}(x,y) = \sup \Big\{u(x) - u(y) \ | \ u \in D(\E^{\rho})_{loc} \cap C(\R^d), \ \Gamma^{\rho}(u,u)   \le   m \ \text{on} \ \R^d  \Big\}.
\]
Here (cf. \eqref{ch5;dfrho}) $\Gamma^{\rho}$ is a symmetric bilinear form on $D(\E^{\rho}) \times D(\E^{\rho})$ such that 
\[
\E^{\rho}(f,g) = \frac{1}{2} \int_{\R^d} \, d\Gamma^{\rho} (f,g).
\]
By \cite[Theorem 4.1]{St5}  and (HP2)
\[
d^{\rho}(x,y) = \| x-y \|, \quad \forall x, y \in \R^d.
\]
Note that by \eqref{ch5;eq;uelliptic} $D(\E^A) = D(\E^{\rho})$.
Fix $x$, $y \in \R^d$. Suppose $(u_n)_{n \ge 1} \subset D(\E^A)_{loc} \cap C(\R^d)$ with $\Gamma^{A}(u_n,u_n)   \le   m \ \text{on} \ \R^d$  is a sequence such that 
\[
d(x,y) = \lim_{n \to \infty} \Big( u_n(x) - u_n(y) \Big). 
\]
Since $  \Gamma^{\rho} \big(u_n / \sqrt{\lambda},u_n / \sqrt{\lambda} \big)  =    \lambda^{-1} \Gamma^{\rho}(u_n,u_n)  \le \Gamma^{A}(u_n,u_n)   \le   m$,
by definition of $d^{\rho}(x,y)$
\[
\| x-y \| = d^{\rho}(x,y) \ge \frac{1}{\sqrt{\lambda}} \lim_{n \to \infty}  \Big( u_n(x) - u_n(y) \Big) =  \frac{1}{\sqrt{\lambda}} d(x,y).
\]
\qed

\begin{remark}\label{ch5;rebmea}
\begin{itemize}
\item[(i)] Assumption (HP2) is only used in order to show \eqref{ch5;eq;intrinsicm} (see proof of Lemma \ref{ch5;l;intrinsic} above).
\item[(ii)] Note that $\tilde{B}_r(x)$ is bounded and open in $\R^d$ for any $x \in \R^d$, $r>0$,  by \eqref{ch5;eq;intrinsicm}.
\end{itemize}
 \end{remark}

The doubling property w.r.t. the intrinsic metric $d(\cdot,\cdot)$ holds: 
\begin{lemma}\label{l;ch5;doublepseu}
Let $n\in \N$ be such that $\lambda\le 2^n$. Then for any $x \in \R^d$, $r>0$
\begin{equation}\label{ch5;eq;doubleintr}
m( \tilde{B}_{2r}(x)  )  \le  C_1^{n+1} \  m(\tilde{B}_r(x) ),
\end{equation}
where $C_1$ is the constant as in \eqref{ch5;doublep}. 
\end{lemma}
\proof
Let $x \in \R^d$, $r>0$.  By Lemma \ref{ch5;l;intrinsic}, it holds
\begin{equation}\label{2.2consequence}
B_{\frac{r}{\sqrt{\lambda}}}(x)\subset  \tilde{B}_{r}(x)  \subset B_{\sqrt{\lambda}r}(x).
\end{equation}
Hence by \eqref{ch5;doublep} and \eqref{2.2consequence}
\[
m( \tilde{B}_{2r}(x)  ) \le m(B_{2 \sqrt{\lambda} r}(x))  \le  C_1 \  m(B_{\sqrt{\lambda} r}(x)).
\]
Using \eqref{ch5;doublep} repeatedly, then \eqref{2.2consequence} and finally the assumption, we get
\[
m(B_{\sqrt{\lambda} r}(x))  \le C_1^{n} \ m(B_{\frac{\sqrt{\lambda}r}{2^n}}(x))\le C_1^{n} \ m(\tilde{B}_{\frac{\lambda r}{2^n}}(x))\le C_1^{n} \ m(\tilde{B}_{r}(x)). 
\]
\qed
\begin{lemma}\label{ch5;l;poincom}
Let $A_1$, $A_2 \in \mathcal{B}(\R^d)$ be bounded sets. Then for any $u \in C(\R^d)$
\[
\int_{A_1} |u - u_{A_1} |^2 \, dm \le \int_{A_1} |u - u_{A_2} |^2 \, dm,
\]
where $u_{B} = \frac{1}{m(B)}\int_{B} u \, dm$, $B \in \mathcal{B}(\R^d)$.
In particular for $A_1 \subset A_2$, we get
\[
\int_{A_1} |u - u_{A_1} |^2 \, dm \le \int_{A_2} |u - u_{A_2} |^2 \, dm.
\]

\end{lemma}
\proof
Let $u \in C(\R^d)$. Then
\begin{eqnarray*}
&& \int_{A_1} |u - u_{A_2} |^2 \, dm - \int_{A_1} |u - u_{A_1} |^2 \, dm \\
&=& -2 u_{A_2}  \int_{A_1} u \, dm +   u_{A_2}^2 \ m(A_1)
+2 u_{A_1}  \int_{A_1} u \, dm -   u_{A_1}^2 \ m(A_1) \\
&=&  -2 u_{A_2} u_{A_1} \ m(A_1)     +   u_{A_2}^2 \ m(A_1)
+2 u_{A_1}^2 \ m(A_1) -   u_{A_1}^2 \ m(A_1) \\
&=& (  u_{A_2} -  u_{A_1}   )^2 \  m(A_1) \ge 0.
\end{eqnarray*}
\qed

Next, we want to show that the scaled strong Poincar\'{e} inequality holds with respect to the intrinsic metric $d(\cdot,\cdot)$. It will be concluded from the next three lemmas in Remark \ref{conclusiverem} below. 
\begin{lemma}\label{ch5;l;poincarebis}
For $x \in \R^d$, $r>0$
\begin{equation}\label{wpoieq}
\int_{\tilde{B}_r(x)} |u - \tilde{u}_{x,r} |^2 \, dm \le c \ r^2 \int_{\tilde{B}_{\lambda r}(x)} \, d\Gamma^{A}(u,u), \quad u \in D(\E^A),
\end{equation}
where $c>0$ is a constant and  $\tilde{u}_{x,r} = \frac{1}{m(\tilde{B}_r(x))}\int_{\tilde{B}_r(x)} u \, dm$.
\end{lemma}
\proof
By (IV) and \eqref{ch5;eq;uelliptic}
for $x \in \R^d$, $r>0$
\begin{eqnarray}
\int_{B_{r} (x)} | u - u_{x,r} |^2  \, dm   &\le& C_3 r^2 \int_{B_{r}(x)}   \|\nabla u \|^2  \, dm,  \notag  \\
 &\le& \lambda \ C_3 r^2 \int_{B_{r}(x)} \, d \Gamma^{A}(u,u), \quad \forall u \in C^{\infty}(\R^d),  \label{ch5;eq; poineuc}  
 \end{eqnarray}
where $C_3$ is the constant as in (IV) and $\lambda$ is the constant of (HP1). Therefore by Lemma \ref{ch5;l;poincom}, \eqref{ch5;eq;intrinsicm}, and \eqref{ch5;eq; poineuc} 
\begin{eqnarray*}
\int_{\tilde{B}_r(x)} |u - \tilde{u}_{x,r} |^2 \, dm &\le&
\int_{B_{\sqrt{\lambda} r} (x)} | u - u_{x, \sqrt{\lambda} r} |^2  \, dm   \le \lambda^2 \ C_3 r^2 \int_{B_{\sqrt{\lambda} r}(x)} \, d \Gamma^{A}(u,u) \\
&\le& \lambda^2 \ C_3 r^2 \int_{\tilde{B}_{\lambda r}(x)} \, d \Gamma^{A}(u,u), \quad \forall u \in C^{\infty}(\R^d).
\end{eqnarray*}
\qed

\begin{lemma}\label{l;2.7jer}
Suppose  $\tilde{c}_1$ and $\tilde{c}_2$ are positive constants such that $\tilde{B}_{r_1}(x_1)$ and $\tilde{B}_{r_2}(x_2)$, $x_1$, $x_2 \in \R^d$ satisfy $\tilde{c}_1 m(\tilde{B}_{r_2}(x_2)) \le m(\tilde{B}_{r_1}(x_1) \cap \tilde{B}_{r_2}(x_2)) \le \tilde{c}_2 m(\tilde{B}_{r_2}(x_2))$ and $u$ is such that  $\int_{\tilde{B}_{r_i}(x_i)} |u - \tilde{u}_{x_i,r_i} |^2 \, dm \le A$, $A>0$, $i=1,2$. Then there exists $\tilde{c}_3=\tilde{c}_3(\tilde{c}_1, \tilde{c}_2)$ such that $\int_{\tilde{B}_{r_1}(x_1) \cup \tilde{B}_{r_2}(x_2)} |u - \tilde{u}_{x_1,r_1} |^2 \, dm \le \tilde{c}_3 A$.
\end{lemma}
\proof
The proof is the same as \cite[Remark 5.4]{Je} with Lebesgue measure replaced by $m$. 
\qed

\begin{lemma}\label{ch5;l;poincare}
The inequality \eqref{wpoieq} implies the scaled weak Poincar\'{e}  inequality:
for $x \in \R^d$, $r>0$
\[
\int_{\tilde{B}_r(x)} |u - \tilde{u}_{x,r} |^2 \, dm \le C \ r^2 \int_{\tilde{B}_{2r}(x)} \, d\Gamma^{A}(u,u), \quad u \in D(\E^A),
\]
where $C>0$ is some constant.
\end{lemma}
\proof
Since the statement  trivially holds if $\lambda \le 2$, we only consider the case $\lambda > 2$. Fix $x \in \R^d$ and $r>0$. For $0 < \alpha \le \frac{\lambda-\varepsilon -1 }{\lambda +1 }$ with small $0 < \varepsilon< 1/\lambda$ depending only on $\lambda$, we can find finitely many points $x_i  \in \overline{\tilde{B}_{(1+\alpha)r}(x)} $, $i=1,\dots,N$, such that $\tilde{B}_{(\alpha + \varepsilon) r} (x_i) \cap \tilde{B}_{r} (x) \ne \emptyset$, $\tilde{B}_{(1+\alpha)r}(x)  \subset  \bigcup_{i=1}^{N} \tilde{B}_{(\alpha + \varepsilon) r} (x_i)$ and $\tilde{B}_{\lambda (\alpha + \varepsilon) r} (x_i) \subset \tilde{B}_{\lambda r} (x)$. Note that using \eqref{ch5;eq;intrinsicm} one can choose the constant $N$ independently of $x$, $r$. Then  by Lemma \ref{ch5;l;poincarebis}, the inclusion $\tilde{B}_{\lambda (\alpha + \varepsilon) r} (x_i) \subset \tilde{B}_{\lambda r} (x)$ and  $\alpha + \varepsilon < 1$, we obtain for each $i = 1,\dots,N$ and $u \in D(\E^A)$
\begin{equation}\label{e;pinlo}
\int_{\tilde{B}_{(\alpha + \varepsilon) r} (x_i)}  |u - \tilde{u}_{x_i ,(\alpha + \varepsilon) r} |^2 \, dm \le c  (\alpha + \varepsilon)^2 \ r^2 \int_{\tilde{B}_{\lambda r}(x)} \, d\Gamma^{A}(u,u) \le c\ r^2 \int_{\tilde{B}_{\lambda r}(x)} \, d\Gamma^{A}(u,u).
\end{equation}
 Applying Lemma \ref {ch5;l;poincom} and then Lemma \ref{l;2.7jer}, \eqref{e;pinlo}, we obtain for $u \in D(\E^A)$
\begin{eqnarray*}
&&\int_{{\tilde{B}_{(1+\alpha)r}(x)}}  |u - \tilde{u}_{x,(1 + \alpha ) r} |^2 \, dm \le
\int_{{\tilde{B}_{(1+\alpha)r}(x)}}  |u - \tilde{u}_{x, r} |^2 \, dm \\
&\le&   \sum_{i=1}^N \int_{{\tilde{B}_{r}(x)}  \cup \tilde{B}_{(\alpha + \varepsilon) r} (x_i) }  |u - \tilde{u}_{x, r} |^2 \, dm  \le
 \tilde{C} N c \ r^2 \int_{\tilde{B}_{\lambda r}(x)} \, d\Gamma^{A}(u,u), 
 \end{eqnarray*}
where $\tilde{C} = \tilde{c}_3$ is the constant of Lemma \ref{l;2.7jer}. Note that $\tilde{C}$ only depends on $\lambda$ and not on $r$ and $x$.  Iterating this argument only finitely many times (depending only on $\lambda$), the statement follows.
\qed

\begin{remark}\label{conclusiverem}
By \cite[Theorem 2.4]{St3} and Lemmas \ref{ch5;l;intrinsic}, \ref{l;ch5;doublepseu}, \ref{ch5;l;poincare}, the scaled strong Poincar\'{e} inequality holds, i.e. for $x \in \R^d$, $r>0$
\[
\int_{\tilde{B}_r(x)} |u - \tilde{u}_{x,r} |^2 \, dm \le c^{\star} \ r^2 \int_{\tilde{B}_{r}(x)} \, d\Gamma^{A}(u,u), \quad u \in D(\E^A),
\]
where $c^{\star}>0$ is some constant.
\end{remark}

\begin{thm}\label{t;conserv}
The Dirichlet form $(\E^A,D(\E^A))$ is conservative.
\end{thm}
\proof
By the doubling property \eqref{ch5;doublep} and \cite[Proposition 5.1, 5.2]{GrHu}
\[
c_{1} \, r^{\alpha^{\prime}} \le m(B_{r}(0)) \le c_{2} \, r^{\alpha}, \quad \forall r \ge 1,
\]
where $c_{1}, c_{2}, \alpha, \alpha^{\prime} > 0$ are some constants. Therefore,
\[
\int_{1}^{\infty} \frac{r}{\log \big (m (B_{r}(0))  \big )} \, dr =\infty.
\]
Then by \cite[Theorem 3.6]{St5} and Lemma \ref{ch5;l;intrinsic} the conservativeness follows  (cf. \cite[Proposition 2.4]{ShTr13b}).
\qed

By Lemmas \ref{ch5;l;intrinsic}, \ref{l;ch5;doublepseu}, \ref{ch5;l;poincare}, $(\E^A,D(\E^A))$ is strongly regular and the properties (Ia)-(Ic) of \cite{St3} are satisfied on $\R^d$. Therefore by \cite[p. 286 A)]{St3} there exists a jointly continuous transition kernel density  $p_{t}(x,y)$ such that
\[
P_t f(x) := \int_{\R^d} p_t(x,y)\, f(y) \ m(dy), \ \  t>0, \ x,y \in \R^d, \ f\in \mathcal{B}_b(\R^d)
\]
is an $m$-version of $T_t f$ if $f  \in  L^2(\R^d , m)_b$. Throughout this paper we set $P_0 : = id$.
Taking the Laplace transform of $p_{\cdot}(x, y)$, we obtain a $\mathcal{B}(\R^d) \times \mathcal{B}(\R^d)$ measurable non-negative resolvent kernel density $r_{\alpha}(x,y)$ such that
\[
R_{\alpha} f(x) := \int_{\R^d} r_{\alpha}(x,y)\, f(y) \, m(dy) \,, \; \alpha>0, \; x \in \R^d, f \in \mathcal{B}_b(\R^d),
\]
is an $m$-version of $G_{\alpha} f$ if $f  \in  L^2(\R^d, m)_b$. For a signed Radon measure $\mu$ on $\R^d$, let us define
\begin{equation*}
R_{\alpha} \mu (x) = \int_{\R^d} r_{\alpha}(x,y) \, \mu(dy) \, , \;\; \alpha>0, \;\;x \in \R^d,
\end{equation*}
whenever this makes sense.

\begin{thm}\label{ch5;l;hkee}
For $x,y \in \R^d$, $t>0$ and any $\varepsilon>0$
\begin{equation}\label{ch5;eq;hkes}
p_t(x,y) \le c \ m (B_{ \sqrt{t}}(x) )^{-1/2} \ m (B_{\sqrt{ t}}(y) )^{-1/2} \ \exp\left(- \frac{\| x-y \|^2}{ \lambda (4+ \varepsilon) t} \right),
\end{equation}
where $c$ is some constant and $\lambda$ is the constant of (HP1).
\end{thm}
\proof
It follows from \cite[Corollary 4.2 and Remarks (ii) in p.286]{St3} that  for $x, y \in \R^d$, $t > 0$ and any $\varepsilon >0$
\begin{equation*}
p_t(x,y) \le c_1\ m \big( \tilde{B}_{\sqrt{t}} (x) \big)^{-1/2} m \big( \tilde{B}_{\sqrt{t}} (y) \big)^{-1/2} \ \exp{ \left ( -\frac{d(x,y)^2}{ (4+\varepsilon)t} \right)},
\end{equation*}
where $c_1$ is some constant. By \eqref{ch5;doublep} and Lemma \ref{ch5;l;intrinsic} the assertion then follows.
\qed

\begin{proposition}\label{ch5;p;strongf}
It holds:
\begin{itemize}
\item[(i)] $(P_t)_{t \ge 0}$ (resp. $(R_{\alpha})_{\alpha > 0}$) is strong Feller, i.e. for $t>0$, we have $P_t(\mathcal{B}_b(\R^d)) \subset C_b(\R^d)$ (resp. for $\alpha>0$, we have $R_{\alpha}(\mathcal{B}_b(\R^d)) \subset C_b(\R^d)$).
\item[(ii)] $P_t (L^1(\R^d,m)_0) \subset C_{\infty} (\R^d)$.
\end{itemize}
\end{proposition}
\proof
Using the transition density estimate \eqref{ch5;eq;hkes}, the statements follow exactly as in \cite[Proposition 3.3]{ShTr13a}.
\qed  

In order to introduce conveniently some notations, we suppose up to the end of this section that there exists a Hunt process 
\begin{equation}\label{ch5;hphp}
\bM : = (\Omega, \F, (\F_t)_{t \ge 0}, \zeta, (X_t)_{t \ge 0}, (\P_x)_{x \in \R^d \cup \{\Delta\} } )
\end{equation}
with transition function $(P_t)_{t \ge 0}$ where $\Delta$ is the cemetery point and the lifetime $\zeta : = \inf \{t \ge 0 \ | \  X_t \in \{\Delta\} \}$.

\begin{remark}
\begin{itemize}\label{ch5;remsuphu}
\item[(i)] The existence of $\bM$ as in \eqref{ch5;hphp} is non trivial. It will be realized in concrete cases via the classical Feller theory through Lemma \ref{ch5;l;fellsem} in Section \ref{ch5;sect3} and via the so-called Dirichlet form method (cf. \cite[Section 4]{AKR},  \cite[Section 2]{ShTr13a}) in Section \ref{ch5;sect4} (see $(\bf{H1})$ and ${(\textbf{H2})^{\prime}}$ there).
\item[(ii)] Under the existence of a Hunt process \eqref{ch5;hphp}, by Theorem \ref{t;conserv} and Proposition \ref{ch5;p;strongf}(i), $\P_x(\zeta = \infty) = 1$ for all $x \in \R^d$.
\end{itemize}
\end{remark}

Let $D$ be an open set in $\R^d$. Then  the part Dirichlet form $(\E^{A,D},D(\E^{A,D}))$ of $(\E^A, D(\E^A))$ on $D$ is a regular Dirichlet form on $L^2(D, m)$ (cf. \cite[Section 4.4]{FOT}). Let  $(T^{D}_t)_{t > 0}$ and $(G^{D}_{\alpha})_{\alpha > 0}$ be the $L^2(D, m)$-semigroup and resolvent associated to the part Dirichlet form $(\E^{A,D},D(\E^{A,D}))$. Let further 
$\sigma_{D} := \inf\{t>0  \ | \ X_t \in D\}$, $P^{D}_{t}f(x) : = \EE_x [f( X_t) ; t<\sigma_{D^c} ]$, $R^{D}_{\alpha}f(x) : = \EE_x \Big[\int_{0}^{\sigma_{D^c}} e^{-\alpha s} f(X_s) \, ds \Big]$, $f \in \mathcal{B}_{b}(D)$. Then $P^{D}_{t} f$ (resp. $R^{D}_{\alpha}f$) is an $m$-version of $T^{D}_t f$ (resp. $G^{D}_{\alpha}f$) for any $f \in L^2(D,m)_b$. Since $P_t^{D} 1_{A}(x) \le P_t 1_{A}(x)$ for any $ A\in \mathcal{B}(D)$, $x \in D$ and $m$ has full support on $\R^d$, $A \mapsto P_t^{D} 1_{A}(x), \; A \in \mathcal{B}(D)$ is absolutely continuous with respect to $1_{D} \cdot m$. Hence there exists a (measurable) transition kernel density $p^{D}_{t}(x, y)$, $x,y \in D$, such that
\begin{equation*}
P_t^{D} f(x) = \int_{D} p_t^{D} (x,y) \,  f(y) \, m(dy) ,\; t>0 \;, \;\; x\in D
\end{equation*}
for $f  \in  \mathcal{B}_{b}(D)$. Correspondingly, there exists a (measurable) resolvent kernel density $r_{\alpha}^{D}(x,y)$, such that
\[
R_{\alpha}^{D} f(x) = \int_{D} r_{\alpha}^{D} (x,y) \,  f(y) \, m(dy) \,, \;\; \alpha>0, \;\;x \in D
\]
for $f  \in \mathcal{B}_{b}(D)$. For a signed Radon measure $\mu$ on $D$, let us define
\[
R_{\alpha}^{D} \mu (x) = \int_{D} r_{\alpha}^{D} (x,y) \, \mu(dy) \, , \;\; \alpha>0, \;\;x \in D
\]
whenever this makes sense. We define $D_A : = \inf \{t \ge 0 \ | \ X_t \in A  \}$ for any $ A\in \mathcal{B}(\R^d)$. Let 
for $\omega \in \Omega$
\begin{equation}\label{ch5;eq;partp}
X_t^{D}(\omega) : = 
\begin{cases}
X_t(\omega) \quad 0 \le t < D_{D^c}(\omega) \\
\Delta \quad t \ge D_{D^c}(\omega).
\end{cases}
\end{equation}
Then $\bM^D : = (\Omega, \F, (\F_t)_{t \ge 0}, (X_t^{D})_{t \ge 0}, (\P_x)_{x \in D \cup \{\Delta\} } )$ is again a Hunt Process by \cite[Theorem A.2.10]{FOT} and its lifetime is $\zeta^D : = \zeta \wedge D_{D^c}$. $\bM^D$ is called the part process of $\bM$ on $D$ and it is associated with the part $(\E^{A,D},D(\E^{A,D}))$ of $(\E^A,D(\E^A))$ on $D$. A positive Radon measure $\mu$ on $D$ is said to be of finite energy integral if
\[
\int_{D} |f(x)|\, \mu (dx) \leq C \sqrt{\E^{A,D}_1(f,f)}, \; f\in D(\E^{A,D}) \cap C_0(D),
\]
where $C$ is some constant independent of $f$ and $\E^{A,D}_1(u,v):=\E^{A,D}(u,v)+\int_D uv\,dm$.
A positive Radon measure $\mu$ on $D$ is of finite energy integral (on $D$) if and only if there exists a unique function $U_{1}^{D} \, \mu\in D(\E^{A,D} )$ such that
\[
\E^{A,D}_{1}(U_{1}^{D} \, \mu, f) = \int_{D} f(x) \, \mu(dx),
\]
for all $f \in D(\E^{A,D}) \cap C_0(D)$. $U_{1}^{D} \, \mu$ is called $1$-potential of $\mu$.  In particular, $R_{1}^{D} \mu$ is a version of $U_{1}^{D} \mu$ (see e.g. \cite[Exercise 4.2.2]{FOT}). The measures of finite energy integral on $D$ are denoted by $S_0^{D}$. We further define $S_{00}^{D} : = \{\mu\in S_0^{D} \, | \; \mu(D)<\infty, \|U_{1}^{D} \mu\|_{\infty, D}<\infty \}$, where  $||f||_{\infty,D}:= \inf \Big\{c>0 \,|\; \int_{D}  1_{\{ \,|f|>c \, \} } \, dm = 0 \Big\}$. If $D = \R^d$, we omit the superscript $D$ and simply write $U_{1}, S_0, S_{00}$.
\begin{proposition}\label{ch5;p;strongfr}
Let $\mu$ be a positive Radon measure and $G \subset \R^d$ relatively compact open. Suppose
\[
\int_G r_1(\cdot,y) \ \mu(dy) \le r_1^G
\]
$\mu$-a.e. on $G$ and $m$-a.e. on $\R^d$, where $r_1^G$ is a continuous function on $\R^d$. Then $1_G \cdot \mu \in S_{00}$.
\end{proposition}
\proof
Since $R_1(1_G \cdot \mu)(x) = \int_G r_1(x,y) \ \mu(dy) \le r_1^G (x)$ for $\mu$-a.e. $x \in G$, $R_1(1_G \cdot \mu) \in L^1(G, \mu)$. This implies that  $1_G \cdot \mu \in S_0$ by \cite[Example 4.2.2]{FOT}. Then $1_G \cdot \mu \in S_{00}$ by \cite[Proposition 2.13]{ShTr13a}.
\qed

\section{Concrete 2-admissible weights with polynomial growth}\label{ch5;sect3}
\begin{defn}
\begin{itemize}
\item[(i)] A function $\psi\in \mathcal{B}(\R^d)$ with $\psi>0$ $dx$-a.e. is said to be a Muckenhoupt $A_2$-weight (in notation $\psi\in A_2$), if there exists a positive constant $A$ such that, for every ball $B\subset\R^d$,
\[
\left(\int_B \psi dx\right) \left(\int_B \psi^{-1} dx\right) \le A \, \left(  \int_{B} 1\, dx \right)^2.
\]
\item[(ii)] A mapping $F: \R^d \to \R^d$ is said to be quasi-conformal if $F$ is one-to-one, the components $F_i$, $i=1,\cdots,d$ of $F$ have distributional derivatives belonging to $L^d_{loc}(\R^d,dx)$ and there is a constant $M>0$ called dilation constant of $F$, so that $dx$-a.e.
\[
\left( \sum_{1 \le i,j \le d} \Big( \partial_{j} F_i \Big)^2 \right)^{1/2} \le M   \left| \emph{det} \ F^{\prime} \right|^{1/d},
\]
where $F^{\prime}(x) = \left( \partial_{j} F_i (x)   \right)_{1 \le i,j \le d}$.
\end{itemize}
\end{defn}
$2$-admissible weights arise typically as in the following example:
\begin{example}\label{ch5;exam2ad} (cf. \cite[Chapter 15]{HKM})
\begin{itemize}
\item[(i)] If $\rho \in A_2$, then $\rho$ is a $2$-admissible weight.
\item[(ii)] If $\rho(x) = | \emph{det} F^{\prime} (x) |^{1-2/d}$ where $F$ is a quasi-conformal mapping in $\R^d$, then $\rho$ is a $2$-admissible weight. 
\end{itemize}
\end{example}

In this section we consider 
 \begin{equation}\label{ch5;eq;2admie}
 \rho(x) = \| x\|^{\alpha}, \quad  \alpha \in (-d, \infty), \quad d \ge 2.
 \end{equation}
Note that for any $\alpha \in (-d, \infty)$ $ \rho(x) = \| x\|^{\alpha}$ satisfies (HP2). In particular, if $\alpha \in (-d,d)$, then $\rho \in A_2$ (see \cite[Example 1.2.5]{Tu}) and if $\alpha \in (-d+2, \infty)$, then $\rho = | \text{det} F^{\prime} |^{1-2/d}$ for some quasi-conformal mapping $F$ (cf. \cite[Section 3]{FKS}). Thus $\rho$ as in \eqref{ch5;eq;2admie} is a 2-admissible weight by Example \ref{ch5;exam2ad}.

\begin{remark}\label{ch5;re2con}
\begin{itemize}
\item[(i)]  The heat kernel estimate \eqref{ch5;eq;hkes} is not explicit. It depends on the volume growth of $m$, hence on $\rho$. In this section, we use the concrete form \eqref{ch5;eq;2admie} for $\rho$ to show the existence of a Hunt process with transition function $(P_t)_{t \ge 0}$ as in \eqref{ch5;hphp} and to find estimates for the 1-potentials corresponding to the drifts of the associated SDE via resolvent kernel estimates. Of course this can be generalized. For instance as in (ii) or by just assuming that the resolvent kernel estimate that we obtain in Lemma \ref{ch5;l;resoes} below is verified for $\rho$.
\item[(ii)]  Let $\phi: \R^d \to \R$ be a measurable function such that  $c^{-1} \le \phi(x) \le c$ $dx$-a.e. for some constant $c \ge 1$. Then by  verifying (I)-(IV), we see that
$\phi \rho$ is a 2-admissible weight, if $\rho$ is a 2-admissible weight. Moreover choosing $\tilde{A} = (\tilde{a}_{ij})$ satisfying (HP1) for $\rho \equiv 1$ we see that $A: = \phi \rho \tilde{A}$ satisfies \eqref{ch5;eq;uelliptic} with respect to the 2-admissible weight $\phi \rho$. In particular, this framework includes Dirichlet forms given as the closure of 
\[
\frac{1}{2} \int_{\R^d} \langle \tilde{A} \nabla f, \nabla g  \rangle \ \phi \rho \, dx, \quad f,g \in C_0^{\infty}(\R^d)
\]
on $L^2(\R^d, \phi \rho dx)$.
\end{itemize}
\end{remark}

\begin{lemma}\label{ch5;l;fellsem}
Let  $\rho$ be as in \eqref{ch5;eq;2admie}.
Then:
\begin{itemize}
\item[(i)] $\lim_{t \to 0} P_t f(x) = f(x)$ for each $x \in \R^d$ and $f \in C_0(\R^d)$.
\item[(ii)] $P_t C_0(\R^d) \subset C_{\infty}(\R^d)$ for each $t>0$.
\end{itemize} 
In particular, $(P_t)_{t \ge 0}$ is a Feller semigroup.  
\end{lemma}
\proof
By Proposition \ref{ch5;p;strongf} (ii), $P_t C_0(\R^d) \subset C_{\infty}(\R^d)$ for each $t>0$. Note that for $\alpha \in [0,\infty)$ and $0< \sqrt{t} \le \|x\|$, we have 
\begin{eqnarray*}
m(B_{\sqrt{t}}(x)) \ge c_d \ (\| x\| - \sqrt{t})^{\alpha}  \sqrt{t}^d,
\end{eqnarray*}
with $c_d = vol(B_1(0))$. Then the statement (i) can be derived as in the proof of \cite[Lemma 3.6 (i)]{ShTr13a}, using Theorem \ref{ch5;l;hkee}, the conservativeness of $(\E^A,D(\E^A))$ and the transition density estimate \eqref{ch5;eq;hkes}. For $\alpha\in (-d,0)$, the proof is the same as in \cite[Lemma 3.6(i)]{ShTr13a}. Hence by \cite[Lemma 2.3]{ShTr13a}  $(P_t)_{t \ge 0}$ is a Feller semigroup.  
\qed

In view of Lemma \ref{ch5;l;fellsem} and the classical Feller theory, there exists a Hunt process 
\begin{equation*}
\bM = (\Omega , \mathcal{F}, (\mathcal{F}_t)_{t\geq0}, \zeta ,(X_t)_{t\geq0} , (\P_x)_{x \in \R^d_{\Delta} }),
\end{equation*}
with state space $\R^d$ and lifetime $\zeta$ such that $P_t(x,B) : = P_t 1_B (x) = \P_x(X_t \in B)$ for any $x \in \R^d$, $B \in \mathcal{B}(\R^d)$, $t \ge 0$. Thus the existence of $\bM$ as in \eqref{ch5;hphp} is guaranteed. As usual any function $f : \R^d \rightarrow \R$ is extended to $\{\Delta\}$ by setting $f(\Delta):=0$.

\subsection{Concrete Muckenhoupt $A_2$-weights with polynomial growth}
In this subsection, we consider the case where $\rho$ as in  \eqref{ch5;eq;2admie} belongs to $A_2$. More precisely, we consider
 \begin{equation}\label{ch5;emucpol}
 \rho(x) = \| x\|^{\alpha}, \quad  \alpha \in (-d, 2), \quad  d \ge 3.
 \end{equation}
\begin{lemma}\label{ch5;l;resoes}
Let  $\rho$ be as in \eqref{ch5;emucpol}.
Then 
\begin{equation}\label{ch5;eq;resoes}
r_1(x,y) \le c_1 \left( \Phi(x,y) + \Psi(x,y)1_{\{ \alpha \in (-d,0)\}}  \right),
\end{equation}
where $\Phi(x,y): = \frac{1}{\|x-y\|^{\alpha +d -2}}$, $\Psi(x,y):= \frac{1}{\|x-y\|^{d -2} \|y\|^{\alpha}}$, and $c_1$ is some constant.
\end{lemma}
\proof
Note that for $\alpha \in [0,2)$, $t >0$, and $x \in \R^d$, we have
\begin{equation*}
 c_2  \ \sqrt{t}^{\alpha + d} \le m(B_{\sqrt{t}}(x)) \le c_3  \ \sqrt{t}^{d} \ (\| x \| + \sqrt{t})^{\alpha},
\end{equation*}
where $c_2$, $c_3$ are some constants. Then the assertion follows as in the proof of \cite[Lemma 3.6 (ii)]{ShTr13a} using the transition density estimate \eqref{ch5;eq;hkes}.
\qed

Define 
\[
V_{\eta} g(x) : = \int_{\R^d} \frac{1}{\| x-y \|^{d-\eta}} \ g(y) \, dy, \quad x \in \R^d, \ \eta > 0,
\]
whenever it makes sense.
\begin{lemma}\label{ch5;l;miz}
Let $\eta \in (0,d)$, $0<\eta - \frac{d}{p} <1$ and $g \in L^p(\R^d, dx)$ with  
\[
\int_{\R^d} (1+\|y\|)^{\eta -d} |g(y)| \, dy < \infty.
\]
Then $V_{\eta}g$ is H\"older continuous of order $\eta - \frac{d}{p}$.
\end{lemma}
\proof
See \cite[Chapter 4, Theorem 2.2]{Miz}.
\qed

\begin{lemma}\label{ch5;l;smresol}
Let $\rho$ be as in \eqref{ch5;emucpol} and $G\subset \R^d$ any relatively compact open set, $p \ge 1$. Suppose
\begin{itemize}
\item[(i)] if  $\alpha \in (-d,-d+2]$, $1_G  \cdot f \ \| x\|^{\alpha} \in L^1(\R^d,dx)$ and $1_G  \cdot f  \in L^q(\R^d,dx)$ with $0 < 2- \frac{d}{q} <1$, 
\item[(ii)]  if $\alpha \in (-d+2,0)$, $1_G  \cdot f \ \| x\|^{\alpha} \in L^p(\R^d,dx)$ with $0 < 2- \alpha - \frac{d}{p} <1$ and $1_G \cdot  f \in L^q(\R^d,dx)$ with $0 < 2- \frac{d}{q} <1$,
\item[(iii)] if $\alpha \in [0,2)$, $1_G  \cdot f \ \| x\|^{\alpha} \in L^p(\R^d,dx)$ with $0 < 2- \alpha - \frac{d}{p} <1$.
\end{itemize}
 Then $R_1(1_G \cdot  |f| m)$ is bounded everywhere (hence clearly also bounded $m$-a.e. on $\R^d$ and $R_1(1_G \cdot |f| m) \in L^1(G,|f|m)$) by the continuous function 
$\int_G |f(y)| \ \left(\Phi(\cdot,y) + \Psi(\cdot,y)1_{\{ \alpha \in (-d,0)\}} \right) \ m(dy)$. In particular, $1_G  \cdot  |f| m \in S_{00}$.
\end{lemma}
\proof
By Lemma \ref{ch5;l;resoes} for any $x \in \R^d$
\begin{eqnarray*}
R_1(1_G \cdot |f|m)(x) &\le& c_1 \ \int_G |f(y)| \ \left(\Phi(x,y) + \Psi(x,y)1_{\{ \alpha \in (-d,0)\}} \right) \, m(dy)\\
&=& c_1 \Big( \ V_{2-\alpha} (1_G \cdot |f| \ \|y\|^{\alpha})(x) + V_2(1_G \cdot |f|)(x) 1_{\{ \alpha \in (-d,0)\}}  \Big).
\end{eqnarray*}
where $c_1$ is the constant as in \eqref{ch5;eq;resoes}. If  $\alpha \in (-d,-d+2]$ and (i) holds, then clearly  $V_{2-\alpha} (1_G \cdot |f| \ \|y\|^{\alpha})$ is continuous. Furthermore, $ V_2(1_G \cdot |f|)$ is continuous by Lemma \ref{ch5;l;miz}. Then by Proposition \ref{ch5;p;strongfr}, $1_G \cdot |f| m \in S_{00}$ (cf. \cite[Proposition 2.13]{ShTr13a}). Thus the statement holds in the case of (i).
The rest of the proof follows from  \eqref{ch5;eq;resoes} as in the proof of \cite[Lemma 3.6 (iii)]{ShTr13a}.
\qed

Up to the end of this subsection, we assume that for each $i,j =1 , \dots ,d$
\begin{itemize}
\item[(HP3)] 
\begin{itemize}
\item[(i)] if  $\alpha \in (-d,-d+2]$, $ \frac{\partial_j a_{ij}}{\rho}  \in   L^1_{loc}(\R^d,m) \cap L^q_{loc}(\R^d,dx)$, $0 < 2- \frac{d}{q} <1$,
\item[(ii)] if $\alpha \in (-d+2,0)$, $ \partial_j a_{ij}  \in L^p_{loc}(\R^d,dx)$ with $0 < 2- \alpha - \frac{d}{p} <1$ and $\frac{\partial_j a_{ij}}{\rho}  \in L^q_{loc}(\R^d,dx)$ with $0 < 2- \frac{d}{q} <1$,
\item[(iii)] if $\alpha \in [0,2)$, $\partial_j a_{ij}   \in L^p_{loc}(\R^d,dx)$ with $0 < 2- \alpha - \frac{d}{p} <1$.
\end{itemize}
\end{itemize}

\begin{lemma}\label{ch5;l;smoothn}
Let  $\rho$ be as in \eqref{ch5;emucpol} and $G\subset \R^d$ any relatively compact open set. Assume (HP1) and (HP3). Then for each $i,j =1 , \dots ,d$
\[
1_{G} \cdot \frac{a_{ii}}{\rho}  m \in S_{00}, \quad 1_{G} \cdot \frac{|\partial_j a_{ij} |}{\rho}  m \in S_{00} .
\] 
\end{lemma}
\proof
For any relatively compact open set $G\subset \R^d$ $ 1_{G} \cdot \frac{a_{ii}}{\rho}  m$  and $1_{G} \cdot \frac{|\partial_j a_{ij} |}{\rho}  m$ are positive finite measures on $\R^d$. Furthermore by \eqref{ch5;eq;uelliptic}, $1_{G} \cdot \frac{ a_{ii} }{\rho} \in \mathcal{B}_b(\R^d)$. Therefore, by  Proposition \ref{ch5;p;strongf} (i) $R_1 \Big(1_{G} \cdot   \frac{a_{ii} }{\rho} m \Big) \in C_b(\R^d)$. Consequently, $1_{G} \cdot \frac{a_{ii}}{\rho}  m \in S_{00}$ (see Proposition \ref{ch5;p;strongfr}).  By the assumption (HP3) and Lemma \ref{ch5;l;smresol}, $1_{G} \cdot \frac{|\partial_j a_{ij} |}{\rho}  m \in S_{00}$.
\qed

We will refer to \cite{FOT} till the end, hence some of its standard notations may be adopted below without definition.
Let  $f^i (x):= x_i$, $i=1,\dots, d$, $x \in \R^d$, be the coordinate functions. Then $f^i \in D(\E^A)_{b,loc}$ and for any $g \in C_0^{\infty}(\R^d)$ the following integration by parts formula holds:
\begin{equation}\label{ch5;eq;ibp}
- \E^A(f^i,g)=    \frac{1}{2} \int_{\R^d} \left( \sum_{j=1}^{d}  \frac{\partial_j a_{ij}}{\rho}  \right) g \, dm, \quad 1 \le i \le d.
\end{equation}
\begin{thm}\label{ch5;t;stfudeco}
 Assume (HP1),  \eqref{ch5;emucpol} (which in particular implies (HP2)), and (HP3). Then it holds $\P_x$-a.s. for any $x \in \R^d$, $i=1,\dots,d$
\begin{equation}\label{ch5;sfd2}
X_t^i = x^i + \sum_{j=1}^d \int_0^t \frac{\sigma_{ij}}{\sqrt{\rho}} (X_s) \, dW_s^j +   \frac{1}{2} \int^{t}_{0}   \left(\sum_{j=1}^d  \frac{ \partial_j a_{ij}}{\rho}\right) (X_s) \, ds, \quad t \ge 0,
\end{equation}\label{ch5;eq;sfdec}
where $(\sigma_{ij})_{1 \le i,j \le d} =  \sqrt{A} $ is the positive square root of the matrix $A$, $W = (W^1,\dots,W^d)$ is a standard d-dimensional Brownian motion starting from zero.
\end{thm}
\proof
By Lemma \ref{ch5;l;smoothn}, \eqref{ch5;eq;ibp}, and \cite[Theorem 5.5.5]{FOT} the strict continuous additive functional, locally of zero energy and corresponding to the coordinate function $f^i \in D(\E^A)_{b,loc}$ is given by 
\[
N^{[f^i]}_t =  \frac{1}{2}  \int^{t}_{0}  \left( \sum_{j=1}^d \frac{ \partial_j a_{ij}}{\rho} \right) (X_s) \, ds, \quad t \ge 0, \quad 1 \le i \le d. 
\] 
The energy measure of $f^i$ denoted by $\mu_{\langle f^i \rangle}$ satisfies $\mu_{\langle f^i \rangle} = \frac{a_{ii}}{\rho} m$.
Therefore by Lemma \ref{ch5;l;smoothn} for any relatively compact open set $G\subset \R^d$, $1_{G} \cdot \mu_{\langle f^i \rangle} \in S_{00}$ and so the positive continuous additive functional in the strict sense corresponding to the Reuvz measure $\mu_{\langle f^i \rangle}$ is given by 
\[
\langle M^{[f^i]}  \rangle_t = \int_0^t  \frac{a_{ii}}{\rho}(X_s) \, ds,
\]
where $M^{[f^i]}_t$ is the continuous local martingale additive functional in the strict sense corresponding to $f^i$. Furthermore since the covariation is
\[
\quad \langle M^{[f^i]}, M^{[f^j]}  \rangle_t = \int_0^t  \frac{a_{ij}}{\rho}(X_s) \, ds,
\]
we can construct a d-dimensional Brownian motion $W$ (on a possibly enlarged probability space $(\overline{\Omega}, \overline{\mathcal{F}}, \overline{\P}_x)$ (cf. \cite[Chapter 3, Theorem 4.2]{KS}) that we call again w.l.o.g. $(\Omega,\mathcal{F}, \P_x)$) such that 
\[
M^{[f^i]}_t = \sum_{j=1}^d \int_0^t \frac{\sigma_{ij}}{\sqrt{\rho}} (X_s) \, dW_s^j,
\]
where $(\sigma_{ij})_{1 \le i,j \le d} =  \sqrt{A} $ is the positive square root of the matrix $A$.
Note that the equation \eqref{ch5;sfd2} holds for all $t \ge 0$ because $(\E^A,D(\E^A))$ is conservative (see Remark \ref{ch5;remsuphu} (ii)).
\qed

\subsection{Concrete weights with polynomial growth induced by quasi-conformal mappings}
In this subsection we consider the case where $\rho$ as in \eqref{ch5;eq;2admie} is induced by the Jacobian of a quasi-conformal mapping.
More precisely, we consider
 \begin{equation}\label{ch5;eqmupo2}
 \rho(x) = \| x\|^{\alpha}, \quad  \alpha \in [2, \infty), \quad d \ge 2.
 \end{equation}
Let 
\begin{equation}\label{ch5;eq;defbk}
B_k : = \{x \in \R^d  \ | \ k^{-1} <  \|x \| < k \}, \quad k \ge 1,
\end{equation}
and for any $G \subset \R^d$
\[
C^{\infty}(\overline{G}) : = \{f:\overline{G} \to \R \ | \ \exists g \in C_0^{\infty}(\R^d), g|_{\overline{G}} = f   \}.
\]
According to \eqref{ch5;eq;uelliptic} the closure of 
\[
\E^{A,\overline{B}_k}(f,g) := \frac{1}{2} \int_{B_k}  \langle  A \nabla f, \nabla g \rangle \, dx, \quad f,g \in C^{\infty}(\overline{B}_k),
\]
in $L^2(\overline{B}_k,m) \equiv L^2(B_k,m)$, $k \ge 1$, denoted by $(\E^{A,\overline{B}_k},D(\E^{A,\overline{B}_k}))$, is a regular Dirichlet form on $\overline{B}_k$ and moreover, it holds:
\begin{lemma}\label{ch5;l;nashhke}
Let  $\rho$ be as in \eqref{ch5;eqmupo2}.
\begin{itemize}
\item[(i)] (Nash type inequality) 
\begin{itemize}
\item[(a)] If $d \ge 3$, then for $f \in D(\E^{A,\overline{B}_k})$
\begin{equation*}
\left\|f\right\|_{2,B_k}^{2 + \frac{4}{d}}\leq c_k \left[\E^{A,\overline{B}_k}(f,f) + \left\|f\right\|_{2,B_k}^2 \right]\left\|f\right\|_{1,B_k}^{\frac{4}{d}}.
\end{equation*}
\item[(b)] If $d=2$, then for $f \in D(\E^{A,\overline{B}_k})$ and any $\delta>0$
\begin{equation*}
\left\|f\right\|_{2,B_k}^{2 + \frac{4}{d+\delta}} \le c_k \left[\E^{A,\overline{B}_k}(f,f) + \left\|f\right\|_{2,B_k}^2 \right]\left\|f\right\|_{1,B_k}^{\frac{4}{d+\delta}}.
\end{equation*}
Here $c_k >0 $ is a constant which goes to infinity as $k \rightarrow \infty$ and $\|f\|_{p,B_k}:=(\int_{B_k}|f|^pdm)^{\frac1p}$, $p\ge 1$.
\end{itemize}
\item[(ii)]
We have for $m$-a.e. $x, y \in B_k$
\begin{itemize}
\item[(a)] if $d \ge 3$, then
\[
r^{B_k}_{1} (x,y)  \le c_1 \frac{1}{\|x-y\|^{d-2}}.
\]
\item[(b)] if $d=2$, then for any $\delta>0$
\[
r^{B_k}_{1} (x,y)  \le c_1 \frac{1}{\|x-y\|^{d+\delta-2}}.
\]
Here $c_1>0$ is some constant.
\end{itemize}

\end{itemize}
\end{lemma}
\proof
Since Sobolev's inequality is applicable on each $B_k$, we can follow the proof of \cite[Lemma 5.4]{ShTr13a} and apply  \eqref{ch5;eq;uelliptic} in 
order to derive the Nash type inequalities in (i). Following  the proof of \cite[Proposition 5.5, Corollary 5.6]{ShTr13a} the assertion (ii) follows.
\qed

Up to the end of this subsection, we assume that 
\begin{itemize}
\item[(HP3)$^{\prime}$]   $\partial_j  a_{ij}  \in L_{loc}^{\frac{d}{2} + \varepsilon}(\R^d,dx)$ for some $\varepsilon>0$ and  each $i,j =1 , \dots ,d$.
\end{itemize}
\begin{lemma}\label{ch5;l;smooloc}
Assume (HP1) and (HP3)$^{\prime}$. Let  $\rho$ be as in \eqref{ch5;eqmupo2} and $f \in L^{\frac{d}{2} + \varepsilon} ( B_k , dx)$ for some $\varepsilon > 0$. Then
\[
1_{B_k} \cdot |f| m \in S^{B_k}_{00}.
\]
In particular 
\[
1_{B_k} \cdot \frac{a_{ii}}{\rho}  m  \in S_{00}^{B_k}, \quad 1_{B_k} \cdot   \frac{|\partial_j a_{ij} |}{\rho}  m \in S_{00}^{B_k}.
\] 
\end{lemma}
\proof
 Using Lemma \ref{ch5;l;nashhke} (ii), Lemma \ref{ch5;l;miz}, and (HP3)$^{\prime}$ the proof is similar to the proof of Lemma \ref{ch5;l;smoothn}, so we omit it (cf. \cite[Lemma 5.8]{ShTr13a}).
\qed

The following integration by parts formula holds for the coordinate functions $f^i\in D(\E^{A,B_k})_{b,loc}$, $i=1,\dots,d$  and $g \in C_0^{\infty}(B_k)$: 
\begin{equation}\label{ch5;pibp}
- \E^{A,B_k}(f^i,g)=    \frac{1}{2}  \int_{B_k}  \left(  \sum_{j=1}^{d} \frac{\partial_j a_{ij}}{\rho} \right)  g \, dm.
\end{equation}

\begin{prop}\label{t;c5lsfd3}
Assume (HP1),  \eqref{ch5;eqmupo2}, and  (HP3)$^{\prime}$. Then the process $\bM$ satisfies 
\begin{equation}\label{ch5;eq;pfd}
X_t^i = x^i + \sum_{j=1}^d \int_0^t \frac{\sigma_{ij}}{\sqrt{\rho}} (X_s) \, dW_s^j +   \frac{1}{2} \int^{t}_{0}   \left( \sum_{j=1}^d \frac{ \partial_j a_{ij}}{\rho} \right)(X_s) \, ds, \quad t < D_{B_k^c},
\end{equation}$\P_x$-a.s. for any $x \in B_k$, $i=1,\dots,d$ where $W=(W^1,...,W^d)$ is a standard d-dimensional Brownian motion starting from zero.
\end{prop}
\proof
Note that \eqref{ch5;eqmupo2} implies (HP2). Applying \cite[Theorem 5.5.5]{FOT} to $(\E^{A,B_k}, D(\E^{A,B_k}))$, the assertion then follows from Lemma \ref{ch5;l;smooloc} and \eqref{ch5;eq;partp}, \eqref{ch5;pibp} (see Theorem \ref{ch5;t;stfudeco} for details).
\qed

\begin{lemma}\label{ch5;l;capli}
Let $\alpha \in [-d+2,\infty)$. Then:
\begin{itemize}
\item[(i)] $\emph{Cap}(\{0\})=0$.
\item[(ii)] For all $x \in  \R^d \setminus \{0\}$
\[
\P_x \Big(\lim_{k \rightarrow \infty} D_{B_k^c} = \infty \Big) = \P_x \Big(\lim_{k \rightarrow \infty} \sigma_{B_k^c} = \infty \Big)=1.
\]
\end{itemize} 
\end{lemma}
\proof
(i) By \cite[Example 3.3.2, Lemma 2.2.7 (ii)]{FOT} and \eqref{ch5;eq;uelliptic}, Cap$(\{0\}) = 0$ if $\alpha \in [-d+2,\infty)$.
(ii) follows from (i), \eqref{ch5;eq;defbk}, Theorem \ref{t;conserv}, and \cite[Lemma 5.10]{ShTr13a}. 
\qed

\begin{thm}\label{ch5;t;solex0}
 Assume (HP1), \eqref{ch5;eqmupo2}, and  (HP3)$^{\prime}$. Then the process $\bM$ satisfies \eqref{ch5;sfd2} for all $x \in \R^d \setminus \{0\}$.
\end{thm}
\proof
Using Lemma \ref{ch5;l;capli} (ii) and \eqref{ch5;eq;pfd} the result follows by letting $k \to \infty$. 
\qed

\begin{remark}\label{ch5;reexshtr}
The results of this section include the particular case where $\phi \equiv 1$ in Remark \ref{ch5;re2con} (ii) with
\begin{equation}\label{ch5;aijtic}
a_{ij}(x) = \tilde{a}_{ij}(x) \|x\|^{\alpha}, \quad \alpha \in (-d,\infty), \quad 1 \le i,j \le d.
\end{equation}
This leads hence to an extension of the results of \cite[Section 3.1 and 3.2]{ShTr13a} with $\phi \equiv 1$ there to the $(a_{ij})$-case. In particular, even if $\tilde{a}_{ij} = \delta_{ij}$ (where $\delta_{ij}$ denotes the Kronecker symbol) we obtain partial improvements of results of \cite[Section 3]{ShTr13a}. For instance by our results it is easy to see that in case $\phi \equiv 1$ \cite[Proposition 3.8 (ii)]{ShTr13a} also holds for $\alpha \in [d, \infty)$, $d \ge 2$. Moreover, in view of Remark \ref{ch5;re2con} (ii) and the results of this section, it is also possible to extend the results of \cite[Section 3.1 and 3.2]{ShTr13a} to the $(a_{ij})$-case with discontinuous $\phi$, $(a_{ij})$ as in \eqref{ch5;aijtic} satisfying (HP3), resp. (HP3)$^{\prime}$.
\end{remark}

\section{Muckenhoupt $A_2$-weights with exponential growth}\label{ch5;sect4}
In this section, we do not use a concrete form of the density estimate \eqref{ch5;eq;hkes}. So rather than considering a concrete $\rho$ as in \eqref{ch5;eq;2admie}, we consider weights in a certain subclass of the Muckenhoupt $A_2$-class. Precisely, we assume the following:
\begin{itemize}
\item[(HP4)]
There exists $\phi \in L^{1}_{loc}(\R^d, dx)$ such that for every cube $Q \subset \R^d$, $d \ge 2$
\[
\frac{1}{dx (Q)} \int_Q e^{|\phi(x) - \phi_Q |} \, dx \le c,
\]
where c is a constant independent of the cube $Q$ and $\phi_Q =  \frac{1}{dx(Q)}  \int_Q \phi \, dx $ and
\begin{equation}\label{ch5;eqexpmuck}
\rho(x) =  e^{\phi(x)}, \quad x \in \R^d.
\end{equation}
\end{itemize}
Then by \cite[IV. Corollary 2.18]{GaFr} $\rho \in A_2$. Consequently $\rho$ is 2-admissible by Example \ref{ch5;exam2ad} (i). Moreover, $\rho$ satisfies (HP2) since for $A_2$-weights $\rho$, it holds $\frac{1}{\rho} \in L^1_{loc}(\R^d,dx)$.\\

In \cite[Section 2]{ShTr13a} we considered a symmetric, strongly local, regular Dirichlet form $(\E,D(\E))$ on $L^{2}(E, m)$ with generator $(L,D(L))$ admitting carr\'e du champ, where $E$ is a locally compact separable metric space and $m$ is a positive Radon measure on $(E, \mathcal{B}(E))$ with full support on $E$.\\

There (with the corresponding objects $(T_t)_{t > 0}$, $(P_t)_{t \ge 0}$, $R_1$, etc., related to $(\E,D(\E))$)     we assumed: 
\begin{itemize}
\item[$(\bf{H1})$] There exists a $\mathcal{B}(E) \times \mathcal{B}(E)$ measurable non-negative map $p_{t}(x,y)$ such that
\[
P_t f(x) := \int_{E} p_t(x,y)\, f(y) \, m(dy) \,, \; t>0, \ \ x \in E,  \ \ f \in \mathcal{B}_b(E),
\]
is a (temporally homogeneous) sub-Markovian transition function (see \cite[1.2]{CW}) and an $m$-version of $T_t f$ if $f  \in  L^2(E , m)_b$.
\item[${(\textbf{H2})^{\prime}}$] We can find $\{ u_n \ | \ n \ge 1 \} \subset D(L) \cap C_0(E)$ satisfying:
\begin{itemize}
\item[(i)] For all $\varepsilon \in \Q \cap (0,1)$ and
$y \in D$, where $D$ is any given countable dense set in $E$, there exists $n \in \N$ such that $u_n (z) \ge 1$, for all $z \in \overline{B}_{\frac{\varepsilon}{4}}(y)$ and $u_n \equiv 0$ on $E \setminus B_{\frac{\varepsilon}{2}}(y)$.
\item[(ii)] $R_1\big( [(1 -L) u_n]^+ \big)$, $R_1\big( [(1 -L) u_n]^- \big)$, $R_1 \big( [(1-L_1)u_n^2]^+ \big)$, $R_1 \big( [(1-L_1)u_n^2]^- \big)$ are continuous on $E$ for all $n \ge 1$ where $L_1$ denotes the $L^1(E,m)$-generator of $(\E,D(\E))$.
\item[(iii)] $R_1 C_0(E) \subset C(E)$.
\item[(iv)] For any $f \in C_0(E)$ and $x \in E$, the map $t \mapsto P_t f(x)$ is right-continuous on $(0,\infty)$.
\end{itemize}
\end{itemize}

Under  $(\bf{H1})$ and ${(\textbf{H2})^{\prime}}$ we showed that there exists a Hunt process with $(P_t)_{t \ge 0}$ as transition function (see \cite[Lemma 2.9]{ShTr13a}). We intend to do the same here in our concrete situation, i.e. we will derive conditions on $a_{ij}$ that imply $(\bf{H1})$ and ${(\textbf{H2})^{\prime}}$.\\

We hence assume in this section that:
\begin{itemize}
\item[(HP5)]  For  $i,j=1,\dots,d$
\[
\partial_{j} a_{ij} \in L^{\infty}_{loc}(\R^d,dx) \quad \text{and} \quad \phi \in L^{\infty}_{loc}(\R^d,dx).
\]  
\end{itemize}

Note that by \eqref{ch5;eq;uelliptic}, $ \left| \frac{a_{ij}}{\rho} \right |  \le  \lambda - \frac{1}{\lambda}(1-\delta_{ij})$, $1 \le i,j \le d$  
and by (HP4) and (HP5) $ \frac{\partial_{j} a_{ij}}{\rho} \in L^{\infty}_{loc}(\R^d,dx)$, $i,j=1,\dots,d$. 
Therefore for $f \in C_0^{\infty}(\R^d)$, we have for the generator $(L^A, D(L^A))$ of $(\E^A, D(\E^A))$
\begin{equation}\label{ch5;dgenerator}
f \in D(L^A) \quad \text{and} \quad L^A f =   \frac{1}{2}\sum_{i,j = 1}^{d} \left(  \frac{a_{ij}}{\rho} \  \partial_{ij} f + \frac{ \partial_{j} a_{ij}}{\rho} \ \partial_i f \right)  \in L^{\infty} (\R^d,m)_0. 
\end{equation}

\begin{thm}\label{ch5;ehuntp}
 Assume (HP1), (HP4) and (HP5). Then there exists a Hunt process $\bM$ satisfying the absolute continuity condition.
\end{thm}
\proof
Using the transition density estimate  \eqref{ch5;eq;hkes}, we can see as in  \cite[Proposition 3.3 (ii)]{ShTr13a} that $(\textbf{H1})$ and \textbf{(H2)}$^{\prime}$ (iii), (iv)  hold. Clearly we can find $\{ u_n \ | \ n \ge 1 \} \subset C_0^{\infty}(\R^d) \subset D(L^A)$ such that $(\textbf{H2})^{\prime}$ (i) is satisfied. Furthermore $(\textbf{H2})^{\prime}$ (ii) for $\{u_n \ | \ n \ge 1 \}$ satisfying $(\textbf{H2})^{\prime}$ (i) follows from \eqref{ch5;dgenerator} and Proposition \ref{ch5;p;strongf} (i).
\qed

Let us write for short
\[
D_k : = B_k(0), \quad k \ge 1.
\]
Note that the $\rho$ is bounded below and above on each $D_k$, $k \ge 1$. Then using Nash type inequalities as in Lemma \ref{ch5;l;nashhke} with $\overline{B}_k$ replaced by $\overline{D}_k$, we obtain for $m$-a.e. $x, y \in D_k$ the resolvent density estimates
\begin{equation}\label{ch5;eq;res1}
r^{D_k}_{1} (x,y)  \le c_1 \frac{1}{\|x-y\|^{d-2}}, \quad \text{if} \ d\ge 3,
\end{equation}
and 
\begin{equation}\label{ch5;eq;res2}
r^{D_k}_{1} (x,y)  \le c_1 \frac{1}{\|x-y\|^{d+\delta-2}},  \quad  \text{for any} \ \delta>0 \quad \text{if} \ d= 2,
\end{equation}
where $c_1$ is some constant.
\begin{lemma}
 Assume (HP1), (HP4), (HP5) and $d \ge 2$. Then:
\[
1_{D_k} \cdot \frac{a_{ii}}{\rho}  m  \in S_{00}^{D_k}, \quad 1_{D_k} \cdot   \frac{|\partial_j a_{ij} |}{\rho}  m \in S_{00}^{D_k}.
\] 
\end{lemma}
\proof
 Using the resolvent density estimates \eqref{ch5;eq;res1}, \eqref{ch5;eq;res2} we can show this similarly to the proof of Lemma \ref{ch5;l;smooloc}. 
\qed

Note that the integration by parts formula \eqref{ch5;pibp} of course holds for $B_k$ replaced by $D_k$.
Consequently, following the proof of Theorem \ref{ch5;t;solex0} we obtain:
\begin{thm}\label{ch5;t;exsol}
 Assume (HP1), (HP4), and  (HP5). Then the process $\bM$ in Theorem \ref{ch5;ehuntp} satisfies \eqref{ch5;sfd2}  $\P_x$-a.s. for any $x \in \R^d$.
\end{thm}

\section{Pathwise unique and strong solutions}\label{ch5;sect5}
In this section we consider 
\begin{itemize}
\item[(HP6)] For each $1 \le i,j \le d$,
\begin{itemize}
\item[(i)]   $ \frac{\sigma_{ij}}{\sqrt{\rho}}$ is continuous on $\R^d$.
\item[(ii)]  $\left \| \nabla \left(  \frac{\sigma_{ij}}{\sqrt{\rho}} \right)  \right \| \in L^{2(d+1)}_{loc} (\R^d,dx)$.
\item[(iii)]  $\sum_{k=1}^{d}  \frac{ \partial_k a_{ik}}{\rho} \in L^{2(d+1)}_{loc} (\R^d,dx)$.
\end{itemize}
\end{itemize}

\begin{thm}\label{ch5;t;ssoleae}
Assume that (HP1), \eqref{ch5;emucpol}, (HP3), and (HP6), resp. (HP1), \eqref{ch5;eqmupo2} (HP3)$^{\prime}$, and (HP6), resp. (HP1), (HP4), (HP5), and (HP6) holds. Then the (weak) solution in Theorem \ref{ch5;t;stfudeco}, resp. Theorem \ref{ch5;t;solex0}, resp. Theorem \ref{ch5;t;exsol} is strong and pathwise unique. In particular, it is adapted to the filtration $(\mathcal{F}_t^W)_{t\geq0}$ generated by the Brownian motion $(W_t)_{t\geq0}$ as in \eqref{ch5;sfd2} and its lifetime is infinite.
\end{thm}
\proof
Assume that (HP1), \eqref{ch5;emucpol}, (HP3), and (HP6), or 
(HP1), \eqref{ch5;eqmupo2}, (HP3)$^{\prime}$, and (HP6),  or (HP1), (HP4), (HP5), and (HP6) holds.
By \cite[Theorem 1.1]{Zh}  under (HP1) and (HP6) for given Brownian motion $(W_t)_{t\geq0}$, $x \in \R^d$ as in \eqref{ch5;sfd2} there exists a pathwise unique strong solution to \eqref{ch5;sfd2} up to its explosion time. The remaining conditions make sure that the unique strong solution is associated to $(\E^A,D(\E^A))$ and has thus infinite lifetime by Remark \ref{ch5;remsuphu} (ii). Therefore the (weak) solution in Theorem \ref{ch5;t;stfudeco},  resp. Theorem \ref{ch5;t;solex0}, resp. Theorem \ref{ch5;t;exsol} is strong and pathwise unique.
\qed
\begin{remark}\label{ch5;renonexpn}
Two non-explosion conditions for strong solutions up to lifetime for a certain class of stochastic differential equations are presented in \cite[Theorem 1.1]{Zh}. For the precise conditions, we refer to \cite{Zh}. By Theorem \ref{ch5;t;ssoleae} and its proof, we know that the solution of \eqref{ch5;sfd2} up to its lifetime fits to the frame of \cite[Theorem 1.1]{Zh}. Therefore, the remaining conditions 
\[
\eqref{ch5;emucpol}, (HP3) \quad \text{or} \quad \eqref{ch5;eqmupo2}, (HP3)^{\prime} \quad \text{or} \quad (HP4), (HP5),
\]  
provide additional non-explosion conditions in \cite[Theorem 1.1]{Zh} for solutions of the form \eqref{ch5;sfd2} that satisfy (HP1) and (HP6). 
\end{remark}

\vspace*{2cm}
Jiyong Shin\\
School of Mathematics\\
Korea Institute for Advanced Study\\
85 Hoegiro Dongdaemun-gu, \\
Seoul 02445, South Korea, \\
E-mail: yonshin2@kias.re.kr\\ \\
Gerald Trutnau\\
Department of Mathematical Sciences and \\
Research Institute of Mathematics of Seoul National University,\\
1, Gwanak-Ro, Gwanak-Gu \\
Seoul 08826, South Korea,  \\
E-mail: trutnau@snu.ac.kr

\end{document}